\documentclass{amsart}
\usepackage{amssymb,amsmath,amsthm,amsfonts,amscd,amstext}
\usepackage[mathscr]{eucal}
\usepackage{array}
\usepackage{graphicx}
\newtheorem{lem}{Lemma}[section]
\newtheorem{prop}[lem]{Proposition}
\newtheorem{thm}[lem]{Theorem}
\newtheorem{cor}[lem]{Corollary}
\newtheorem{df}{Definition}[section]
\newtheorem{rem}[lem]{Remark}

\newtheorem{question}{Question}[section]

\begin{document}
\title{A rough analytic relation on partial differential equations}
\author{Tsuyoshi Kato}
\address{Department of Mathematics, Graduate School of Science, Kyoto
University, Sakyo-ku, Kyoto 606--8502, Japan}
\email{tkato@math.kyoto-u.ac.jp}
\author{Satoshi Tsujimoto}
\address{Department of Applied Mathematics and Physics, Graduate School
of Informatics, Kyoto University, Sakyo-ku, Kyoto 606--8501, Japan}
\email{tujimoto@i.kyoto-u.ac.jp}
\begin{abstract}
We introduce some analytic relations on the set of partial differential equations of two variables.
It relies on a new comparison method to give rough asymptotic estimates for
solutions which obey   different partial differential equations. It uses a kind of scale transform called 
tropical geometry which connects automata with real rational dynamics.
Two different solutions can be considered when their defining equations are transformed to the same automata at infinity.
We have a systematic way to construct related pairs of different partial differential equations,
and also construct some unrelated pairs concretely. 
These verify that the new relations are non trivial.
We also make numerical calculations and compare the results for both related and unrelated pairs of PDEs.
\end{abstract}

%
%
%
\maketitle

\section{Introduction}
Scaling limits play one of the central roles in
discrete dynamical systems, which create
another dynamics with different analytic aspects.
Tropical geometry and ultradisctete dynamical systems arose from very different 
contexts each other, and have been developed quite separately  [V], [TTMS].
However from the view point of scaling limits, their method surprisingly coincides, 
which connect real rational dynamics with  piece-wise linear systems.

Such scaling limit provides with an important prototype
in mathematical physics. 
The Korteweg-de Vries (KdV) equation
$u_s +uu_x + u_{3x}=0$ is a PDE which 
has been derived from the discrete Lotka-Volterra equation
$ z_{N+1}^{t+1} = z_{N+1}^{t} (1+  z_N^{t+1} )(1+z_{N+2}^{t})^{-1}$
\cite{hirota,hirotsuji}. 
The  above  scaling limit
changes the discrete dynamics into   the Lotka-Volterra cell automaton
$ U_{N+1}^{t+1} = U_{N+1}^{t} 
 + \max(0, U_N^{t+1}) - \max(0, U_{N+2}^{t})$.
In the light of integrable systems, these three dynamics share
 common features, like  existence of soliton solutions, many conserved
 quantities and so on [TTMS]. 
 
 These scaling limits allow us to study analytic aspects at the same time 
 for  three categories of dynamical systems which sit in different hierarchies mutually.
The classes of dynamics which can be analyzed by these scaling limits
are rather broad beyond integrable systems.
 In fact  it is  easy  to find discrete dynamical systems which are far from  integrable,
 but which are transformed into  integrable cell automata.
 (see section \ref{sec:3}.C). 
If we consider a situation 
when  two discrete dynamical systems are transformed into the same integrable cell automaton,
  where one is integrable and the other is not, then 
  it will be quite natural to try to produce some systematic method to analyze 
  various PDEs which include wide classes of dynamical systems as above.
 
Arising from a new comparison method to 
 study rough asymptotic growth among solutions to different partial differential equations,
in this paper we  introduce analytic relations on the set of
partial differential equations of two variables. As a first step 
 in this paper we show  non triviality of such relations.

Let
$P_1(\epsilon,u,u_x,u_s, \dots)$ and $P_2(\epsilon,v,v_x,v_s, \dots)$ be  two different partial differential equations
for dependent variables  $u(x,s)$ and $v(x,s)$ respectively.
Our basic method  is to compare solutions between $P_1$ and $P_2$
in terms of the initial data and the distorsions of the higher derivatives of their solutions.

For small  parametrization $\epsilon >0$, let: 
\begin{align*}
{\bf PDE}_2 =\{ P(\epsilon, u,u_x,u_s, &  u_{xs}, \dots, u_{\alpha s})=0 : \\
& P  \text{ are polynomials with real coefficients} \}
\end{align*}
be all the set of families of polynomial type partial differential equations with 2 variables 
parametrized by $\epsilon>0$.

In this paper we introduce some analytic relations $\sim$ between elements in ${\bf PDE}_2$.
They are given by {\bf uniform asymptotic estimates} for all positive solutions of different PDEs.
We will focus on two aspects, where one is growth rates of solutions and the other is 
the domains of solutions on the space variable.
For the first we will treat two growth cases,  exponential $\sim^e $ and double-exponential $\sim^{e^e}$.
For the second we also treat two cases, 
where one is  $\sim_{fin}$ which include  domains of bounded space variables, 
and the other is $\sim_{\infty}$ which  restricts only on infinite domains in space.
Thus in total there are four classes of the relations:
$$\begin{matrix} 
& \sim_{fin}^e & \leq &  \sim_{fin}^{e^e} \\
&&& \\
&| \wedge && | \wedge  \\
&&& \\
&   \sim_{\infty}^e & \leq  & \sim_{\infty}^{e^e}
\end{matrix}
$$
where 
 $\sim_*^+ \  \  \leq \ \  \sim_{*'}^{+'} $ means that $u  \ \sim_*^+ \ v$ implies $u \ \sim_{*'}^{+'} \ v$ for any two solutions.

In this paper we introduce two subsets: 
$$ {\bf PDE}_2^{fin}     \subset {\bf PDE}_2^{\infty} \subset {\bf PDE}_2.$$
There are also some stratifications  $(M,c,D,L; \alpha)$ and  $(M,c,L,k,D)$ over both $\sim_*^+$ and
${\bf PDE}_2^*$ with respect to some constants respectively ($2.B, 3.B$).

Our main aim here is to verify non triviality of the relations over these spaces.
In [K2], we have introduced some constructive way to obtain partial differential equations,
which produces many related pairs with respect to the classes:
$$\ ({\bf PDE}_2^{fin} , \sim_{fin}^*) , \quad  \ ({\bf PDE}_2^{\infty} , \sim_{\infty}^*)$$
where $*=$  $e$ or $e^e$.

Let us explain more details.
Our basic idea is to approximate PDE by discrete rational dynamics
 of the form:
$$z_{N+1}^{t+1}= f(z_{N-l_0}^{t+1}, \dots, z_N^{t+1}, z_{N- l_1}^t ,\dots, z_{N+k_1}^t, \dots, 
z_{N+k_{d+1}}^{t-d})$$
by introducing  scaling parameters $z_N^t = \epsilon^m u(x,s)$ and
$(N,t)=( \epsilon^{-p}x, \epsilon^{-q}s)$.
Notice that for any $n$ variable function $f$ as above, one can choose various 
types of the sets of variables $(z_{N-l_0}^{t+1}, \dots, z_N^{t+1}, z_{N- l_1}^t ,\dots, z_{N+k_1}^t, \dots, z_{N+k_{d+1}}^{t-d})$.
We say that the sets:
$$\{(z_{N-l_0}^{t+1}, \dots, z_N^{t+1}, z_{N- l_1}^t ,\dots, z_{N+k_1}^t, \dots, z_{N+k_{d+1}}^{t-d}) ,
(m,p,q)\}$$ are the {\em approximation data}.

Let us say that the above rational dynamics is {\em consistent}, 
if all $k_1, \dots, k_{d+1} \leq 1$ hold.
 As a general procedure, the rational dynamics with the scaling parameters above give 
pairs of partial differential equations 
${\bf F}( \epsilon, u,u_x, \dots)=0$ as the leading terms,
and the error terms ${\bf F}^1(\epsilon, u,u_x, \dots)$
by use of Taylor expansions ($3.B$):
 \begin{align*}
  & z_{N+1}^{t+1}- f(z_{N-l_0}^{t+1}, \dots, z_{N+k_{d+1}}^{t-d}) 
   =  \epsilon^m{\bf F }( \epsilon , u,u_x,  u_s, u_{xs}, \dots , u_{\alpha x}, u_{\alpha s}) \\
& \quad  + \epsilon^{m+1}{\bf F}^1(\epsilon, u,u_x, u_s, u_{xs}, \dots, u_{\alpha s}, 
  u_{(\alpha +1)x}(\xi_{\alpha+1,0}), \dots, 
 u_{(\alpha +1)s}(\xi_{0, \alpha+1})  )
  \end{align*}
  We say that a partial differential equation 
$P(u,u_x,  u_s, u_{xs}, \dots , u_{\alpha x}, u_{\alpha s}) $ is in ${\bf PDE}_2^{\infty}$, if there is 
an induced  pair $({\bf F}, {\bf F}^1)$ as above and a positive function $h >0$ so that:
(1)  ${\bf F}$ satisfies the equality:
$${\bf F}(\epsilon, u,u_x,  u_s, u_{xs},   \dots , u_{\alpha x}, u_{\alpha s}) 
= \frac{P(u,u_x,  u_s, u_{xs}, \dots , u_{\alpha x}, u_{\alpha s})}{h(u,u_x,  u_s, u_{xs}, \dots , u_{\alpha x}, u_{\alpha s})} $$
(2) 
 there is a constant $C\geq 0$ so that the pointwise estimates hold:
\begin{align*}
|{\bf F}^1  (\epsilon, u,u_x, u_s, &   \dots,  u_{\alpha s},  u_{(\alpha +1)x}(\xi_{\alpha+1,0}),  u_{\alpha x,s}(\xi_{\alpha,1}),
\dots,  u_{(\alpha +1)s}(\xi_{0, \alpha+1})  )| \\
 &  \leq C ( |u_{(\alpha +1)x}(\xi_{\alpha+1,0})| +  | u_{\alpha x,s}(\xi_{\alpha,1})| + \dots +
| u_{(\alpha +1)s}(\xi_{0,\alpha+1}) | ).
\end{align*}
 $P$ is said to be in ${\bf PDE}_2^{fin}$, 
if in addition the discrete dynamics is consistent.

  For the domains of solutions,  we consider on:
$$  \begin{cases}
   (0, A_0) \times [0, T_0) & P \in {\bf PDE}_2^{fin},\quad A_0, T_0 \in (0, \infty], \\
(0, \infty) \times [0, T_0) & P \in {\bf PDE}_2^{\infty}, \quad T_0 \in (0, \infty].
\end{cases}$$

The following result  was obtained by the analysis of the combination of Lipschitz geometry with some scaling limits 
called tropical geometry:
\begin{thm}[K2]
Suppose two PDEs $P,Q \in {\bf PDE}_2^*$ are obtained from two
relatively elementary and increasing functions $f$ and $g$ which are mutually  tropically equivalent,
where $*=$ fin or $\infty$.
Then they are mutually related:
$$ P \sim_*^+ Q \qquad \text{ if } P,Q \in {\bf PDE}_2^*$$
where  $+ = e$ or $e^e$, and
  $+ $ can be $ e$ only when both $f$ and $g$ are linear.
\end{thm}

As a concrete case, we have the following:

\begin{prop}  
Two partial differential equations of order $1$:
$$v_s+ \frac{\epsilon}{2} vv_x - \frac{1}{2} v^2=0, \quad 2 u_s  + \frac{\epsilon}{2} u(u_s+ u_x)  =0$$
are both in ${\bf PDE}_2^{\infty}$ in the class $(20,2,1,2,1)$,
and they are mutually related in $\sim_{\infty}^{e^e}$ in the class $(40, 4, 1,2; 1)$.
\end{prop}

The following  theorem suggests that ``sizes of   \ ${\bf PDE}_2^* / \sim_*^+$ \  will be large'':

    \begin{thm}   
(1)  For any $M,D,L \geq 1$, 
  there are  $l_0$ and $a_0,b_0$ so that for each  even $l=2m  \geq l_0$,
  two PDEs $u,v: (0, a_0) \times [0, b_0) \to (0, \infty)$ of order $l$:
  \begin{align*}
&  \epsilon u_s + \frac{\epsilon^2}{2} u_{2s} + \dots + \frac{\epsilon^l}{l!} u_{ls} + \epsilon^2 u_{xs} 
+\frac{\epsilon^3}{2} u_{xss} + \dots + \frac{\epsilon^l}{l!} u_{(l-1)M's}=0, \\
&  \epsilon(v_s+v_x) + \epsilon^2 (v_{xs} + \frac{1}{2}v_{2s} + \frac{1}{2} v_{2x}) + \dots + \frac{\epsilon^l}{l!} v_{lx} =0
\end{align*}
  are both  in ${\bf PDE}_2^{fin}$ in the class $(1, 1, 0, 1, 1)$, and
  they   are mutually unrelated in $\sim_{fin}^e$ in  the  classes $(M,1, D, L; l)$.

\vspace{5mm}

(2)
For any $M, c, D, L \geq1$, there is $I_0$ so that for all $I \geq I_0$,
 pairs of PDEs given  by:
$$ 4u_s +\epsilon u(u_s+u_x)=0,  \quad v_s=I v_x$$
are both in ${\bf PDE}_2^{\infty}$  in $(20, 2,1,I+1, 1)$, and 
are mutually unrelated in $\sim_{\infty}^{e^e}$ in the class $(M, c, D, L;1)$.
\end{thm}
 The analysis for (2) touches with the technique in the field of  the {\em singular perturbations}
 in the sense of continuity of solutions at $\epsilon =0$.

The discrete dynamics used in the proof of (1) above satisfy  the following property (def $4.1$):
\begin{lem}
There is a  pair of the discrete dynamics which  are mutually  infinitely unrelated.
  \end{lem}

Now in particular one has obtained the followings:
\begin{cor}
There are  unrelated pairs both in 
$$(1) \ \ ( {\bf PDE}_2^{fin} ,\ \sim_{fin}^e), \ \
(2)   \ \ ( {\bf PDE}_2^{\infty} ,  \ \sim_{fin}^{e^e})$$ 
with respect to any constants.
\end{cor}

\vspace{3mm}

Finally we would like to address some problems which arise from the line of this paper:
\begin{question}
Are there unrelated pairs in the following classes:
$$(1) \ \ ( {\bf PDE}_2^{fin} ,\ \sim_{fin}^{e^e}), \ \ 
(2)  \ \ ( {\bf PDE}_2^{\infty} ,\ \sim_{\infty}^e), \ \ 
(3) \ \ ( {\bf PDE}_2^{\infty} ,\ \sim_{\infty}^{e^e}).$$
\end{question}

\vspace{3mm}

In order to approach such problems, 
 computer calculations seem very convenient.
In fact in many examples estimated here, one can express very concrete constants
 even though they are far from optimal.
 So in section 5, we have given computer calculation results.
 It turns out that certainly their data reflect related and unrelated situations.

In order to treat  the real models in
physics, we need to extend this method for systems of partial differential equations.
In particular the Korteweg-de Vries (KdV) equation is
 quite intriguing systems and  has surprising relationship with 
the soliton cellular automaton system known as Box and ball system
(BBS) through the time-discretized Lokta-Volterra system.
We would expect that study of such direction might lead us to rough classifications
of non-equilibrium systems.

\section{Classes for partial differential equations}\label{sec:classPDE}
{\bf 2.A Initial conditions:}
Let us consider $C^{\alpha+1}$ functions 
$u, v: (0, A_0) \times [0, T_0) \to (0, \infty)$, where $A_0, T_0 \in (0, \infty]$.

For  $0< \epsilon <1$,
we  introduce the {\em initial rates}:
$$[u:v]_{\epsilon} \equiv \sup_{(x,s)  \in 
(0, A_0) \times [0, \epsilon] \cup (0,  \epsilon] \times [0, T_0)}
 (\frac{u(x,s)}{v(x,s)})^{\pm 1}.$$

Let $u: (0, A_0)  \times[ 0, T_0) \to (0, \infty)$
be a function of   class $C^{\alpha + 1}$ as above.
Then we  introduce the uniform norm  of $u$ of order $\alpha+1$ by:
$$||u||_{\alpha+1} = \max_{ \partial_i = \  \partial_x , \ \partial_s} \quad
\{ ||\frac{\partial^{\alpha+1} u}{\partial_1 \dots \partial_{\alpha+1}}||_{C^0((0, A_0)  \times [0, T_0)) } \}.$$

Suppose $u(x,s) >0$ is positive.
Then we introduce the {\em higher derivative rates} by:
 $$K_{\alpha+1}(u)  \equiv  \frac{||u||_{\alpha+1}}{\inf_{(x,s) \in (0, A_0)  \times [0, T_0)} u(x,s)}.$$

\vspace{3mm}

{\bf 2.B Analytic relations on  partial differential equations:}
Let:
$$(*,+) \in \{(fin, e) , \ (fin , e^e), \ (\infty , e), \ (\infty , e^e)\}.$$

\begin{df} Let $P , Q \in {\bf PDE}_2^*$ of order $\alpha$. 
$P$ and $Q$ are related in $(*,+)$:
$$P \ \sim_*^+ \ Q$$
if there are constants $M,c, D, L \geq1$ and 
 $C,C'$ so that
for any positive solutions 
$u, v : (0, A_0) \times [0, T_0) \to (0, \infty)$:
$$P(\epsilon, u,u_x,u_s, u_{xs}, \dots, u_{\alpha s})=0, \quad
Q(\epsilon, v,v_x,v_s,v_{xs}, \dots, v_{\alpha s})=0$$
they satisfy the asymptotic estimates:
$$ (\frac{u(x,s)}{v(x,s)})^{\pm 1} \leq 
\begin{cases}
M^{c^{\epsilon^{-D}(x+s+1)}}([u:v]_{L \epsilon})^{c^{\epsilon^{-D}(x+s+1)}}  &  + =e^e, \\
M^{\epsilon^{-D}(x+s+1)}[u:v]_{L \epsilon} & +=e
\end{cases}$$
 for all $0 < L \epsilon \leq \min(\frac{1}{C K}, A_0, T_0, C')$,
where $K = \max( K_{\alpha+1}(u), K_{\alpha+1}(v))$.
\end{df}

The following is immediate:
\begin{lem}
Let $P,Q \in {\bf PDE}_2^{\infty}$. 
$P_1\sim_{\infty}^+ P_2$ holds whenever 
$P_1\sim_{fin}^+ P_2$ for both $+= e$ or $e^e$.
\end{lem}

When we specify the constants, then we denote:
$$ P \sim_*^{e^e} Q \quad \text{ in } (M,c,D,L ; \alpha), \ \ \text{ or }
 P \sim_*^e Q \quad  \text{  in } (M,1,D,L; \alpha).$$

The following is immediate:
\begin{lem}
 Suppose $M \leq M'$, $c\leq c'$, $D\leq D'$ and $L \leq L'$  hold. Then 
 $P_1 \sim_*^+ P_2$ in $(M,c,D,L ; \alpha)$ implies $P_1\sim_*^+ P_2$ in $(M',c', D' ,L'; \alpha)$.
 \end{lem}
 
 \vspace{3mm}
 
\begin{rem}Notice that even though two pairs 
 $P_1 , P_2 \in {\bf PDE}_2^*$ are unrelated in $\sim_*^+$  in the class $(M,c,D,L; \alpha)$, 
they might still be related in $\sim_*^+$ in the class $(M',c',D',L' ; \alpha)$.
\end{rem}

\vspace{3mm}

In [K2], we have obtained a general method to produce related partial differential equations.
In the next section we sketch its construction, and
in section 3.D.2, we calculate a  concrete example of a pair of  related PDEs 
which arise from the method.

 \section{Construction of related PDEs}\label{sec:3}
{\bf 3.A Evolutional discrete dynamics:}
Let $f =\frac{k}{h}$ be a rational function of $n$ variables, where $k$ and $h$ are both polynomials.
An evolutional discrete dynamics is given by flows of the form
$\{ z_N^t\}_{t,N \geq 0}$, where one regards that $t$ is time parameter:
$$\begin{matrix}
z_0^0 & z_1^0 & \dots & z_N^0 & \dots &  \downarrow  ( t=0)\\
 z_0^1 & z_1^1 & \dots & z_N^1 & \dots  & \downarrow (t=1) \\
 \dots & \dots & \dots & \dots & \dots & \\
 z_0^t & z_1^t & \dots & z_N^t & \dots & \\
 \dots & \dots & \dots & \dots & \dots &
\end{matrix}$$

A general equation of evolutional discrete dynamics is of the form:
$$z_{N+1}^{t+1}= f(z_{N-l_0}^{t+1}, \dots, z_N^{t+1}, z_{N- l_1}^t ,\dots, z_{N+k_1}^t, \dots, 
z_{N+k_{d+1}}^{t-d})$$
 where $l_i, k_j \geq 0$,  $N \geq \max(l_0, \dots, l_{d+1})$ and $ t \geq d$,
with initial values: 
$$\bar{z}_0^0 \equiv 
\{z_a^t\}_{ 0 \leq a \leq \max(l_0,\dots, l_{d+1}) , t =0,1, \dots} 
 \cup \{z_N^h\}_{0 \leq h \leq d, N=0,1, \dots}.$$

\vspace{3mm}

 {\bf 3.B Approximation of PDE by rational dynamics:}
Let us consider $C^{\alpha+1}$ functions 
$u: (0, A_0 ) \times [0, T_0) \to (0, \infty)$ where $A_0, T_0 \in (0, \infty]$.
We consider its Taylor expansion,  where  $|(x,s)-\xi_{ij}| \leq | ( i \epsilon^p , j \epsilon^q)|$:
 \begin{align*}
z_{N+i}^{t+j} & =   u  (x  + i \epsilon^p   , s+ j \epsilon^q)
   = u+i  \epsilon^p u_x +  j \epsilon^q u_s+
    \frac{(i\epsilon^p)^2}{2} u_{2x} +  \frac{(j \epsilon^q)^2}{2} u_{2s} \\
    & \qquad +j \epsilon^q i \epsilon^p u_{xs} +
    \dots  + \frac{(i\epsilon^p)^{\alpha}}{\alpha !} 
    u_{\alpha x} + \frac{(j \epsilon^q)^{\alpha}}{\alpha !} u_{\alpha s} \\
&  \qquad \quad  + \frac{(i\epsilon^p)^{(\alpha +1)}}{(\alpha +1)!} u_{(\alpha +1)x}(\xi_{ij}) +
  \dots
 + \frac{(j\epsilon^q)^{(\alpha +1)}}{(\alpha +1)!} u_{(\alpha +1)s}(\xi_{ij}) .
 \end{align*}

Let $f = \frac{k}{h}$ be as in $3.A$ and  choose an approximation data.
Then we consider the corresponding  discrete dynamics
$z_{N+1}^{t+1}= f(z_{N-l_0}^{t+1}, \dots,  z_{N+k_{d+1}}^{t-d})$ 
and change of variables 
  $ \epsilon^m u(x,s) = z^t_N$, 
$N = \frac{x}{\epsilon^p} $ and $  t = \frac{s}{\epsilon^q}$.
If $A_0 < \infty$, then we assume $k_i \leq 1$ for all $1 \leq i \leq d+1$.

Let us insert the expansion and take their difference:
\begin{align*}
& z_{N+1}^{t+1}- f(z_{N-l_0}^{t+1}, \dots, z_{N+k_{d+1}}^{t-d}) \\
& =  \epsilon^m u(  x+ \epsilon^p, s+ \epsilon^q) - 
 f( \epsilon^m u(x - l_0 \epsilon^p, s+ \epsilon^q),  \dots, 
 \epsilon^m u(x+ k_{d+1} \epsilon^p, s - d \epsilon^q)) \\
& =  \frac{ \epsilon^m  F^1(u) + \epsilon^{m+p} F^2( u_x) +   \epsilon^{m+q} F^3( u_s) +
\epsilon^{2m+p} F^4( u, u_x)  + \dots  }{
 h( \epsilon^m u(x - l_0 \epsilon^p, s+ \epsilon^q),  \dots, 
 \epsilon u(x+ k_{d+1} \epsilon^p, s - d \epsilon^q ))}  
 \end{align*}
 where $F^i(u,u_x, \dots)$ are monomials.

 Let us divide its numerator into two parts: 
 \begin{align*}
&  \epsilon^m  F^1(u) + \epsilon^{m+p} F^2( u_x) +   \epsilon^{m+q} F^3( u_s) +
\epsilon^{2m+p} F^4( u, u_x)  + \dots   \\
& = \epsilon^m P(\epsilon, u,u_x,  u_s, u_{xs}, \dots ,  u_{\alpha s}) + \\
& \qquad \quad  \epsilon^{m+1}Q(\epsilon, u, .. , \{u_{(\alpha +1)x}(\xi_{\alpha+1,0}),u_{\alpha x ,s}(\xi_{\alpha,1}),
\dots ,  u_{(\alpha +1)s}(\xi_{0, \alpha+1})\} )
\end{align*}
so that the first term contains  only monomials whose derivatives of $u$ are up to order $\alpha$,
and  all the monomials of the second term contain derivatives of $u$ of order $\alpha+1$. 
Then putting:
\begin{align*}
& {\bf F}( \epsilon , u,u_x,  u_s, u_{xs}, \dots , u_{\alpha x}, u_{\alpha s}) \\
& \qquad \qquad = P( \epsilon,  u,u_x,  \dots, u_{\alpha s})/ 
h( \epsilon^m u(x - l_0 \epsilon^p, s+ \epsilon^q),  \dots)  , \\
&  {\bf F}^1( \epsilon , u, .., \{u_{ix,js}(\xi_{ij}) \}_{i+j=\alpha+1})\\
& \quad = Q(\epsilon, u, .. , \{ u_{ix, js}(\xi_{ij})\}_{i+j=\alpha+1} )/ 
h( \epsilon^m u(x - l_0 \epsilon^p, s+ \epsilon^q),  \dots),
\end{align*}
 one obtains the following expression:
  \begin{align*}
   z_{N+1}^{t+1}- & f(z_{N-l_0}^{t+1}, \dots, z_{N+k_{d+1}}^{t-d}) \\
&   = \epsilon^m  {\bf F }( \epsilon , u,u_x,  u_s, u_{xs}, \dots , u_{\alpha x}, u_{\alpha s}) \\
& \qquad  \quad + \epsilon^{m+1}{\bf F}^1(\epsilon, u,u_x, \dots, u_{\alpha s}, \{u_{ix,js}(\xi_{ij}) \}_{i+j=\alpha+1}).
  \end{align*}

\begin{df}  
A partial differential equation 
$P(u,u_x,  u_s, u_{xs}, \dots , u_{\alpha x}, u_{\alpha s}) $ is in ${\bf PDE}_2^{\infty}$, if there is 
an induced  pair $({\bf F}, {\bf F}^1)$ as above and a positive function $h >0$ so that:

(1)  ${\bf F}$ satisfies the equality:
$${\bf F}(\epsilon, u,u_x,  u_s, u_{xs}, \dots , u_{\alpha x}, u_{\alpha s}) 
= \frac{P(\epsilon, u,u_x,  u_s, u_{xs}, \dots , u_{\alpha x}, u_{\alpha s})}{h( \epsilon^m u(x - l_0 \epsilon^p, s+ \epsilon^q),  \dots)} $$

(2) 
 there is a constant $C\geq 0$ so that the pointwise estimates hold:
\begin{align*}
|{\bf F}^1  (\epsilon, u,u_x, u_s, &  \dots, \{u_{ix,js}(\xi_{ij}) \}_{i+j=\alpha+1})| \\
 &  \leq C ( |u_{(\alpha+1)x}(\xi_{\alpha+1,0})| + \dots +| u_{(\alpha +1)s}(\xi_{0, \alpha+1}) | ) .
\end{align*}

 $P$ is in ${\bf PDE}_2^{fin}$
if in addition, the discrete dynamics is consistent.
\end{df}

Let us consider the last term with the estimate  (see $2.A$):
$$C ( |u_{(\alpha+1)x}(\xi_{\alpha+1,0})| + \dots +| u_{(\alpha +1)s}(\xi_{0, \alpha+1}) | )
 \leq  \ Cl \  ||u||_{\alpha+1}$$
where $l$ is the number of the summation of $\alpha+1$ derivatives of $u$.
We say that the number $Cl $ is the {\em error constant} for the approximation of $P$.
Notice that the error constant is determined by the approximation data and the original rational function $f$.

\vspace{3mm}

For $P \in {\bf PDE}_2^{fin}$, we consider $(0, A_0) \times [0, T_0)$ as domains of solutions,
and $(0, \infty) \times [0, T_0)$ for $P \in {\bf PDE}_2^{\infty}$, 
where $A_0, T_0 \in (0, \infty]$.

\vspace{3mm} 

 Let: 
 $$\{(z_{N-l_0}^{t+1}, \dots, z_N^{t+1}, z_{N- l_1}^t ,\dots, z_{N+k_1}^t, \dots, z_{N+k_{d+1}}^{t-d}) ,
(m,p,q)\}$$ be the approximation data. Then one obtains several numbers:
$$L =\max(l,d), \ \ D=\max(p,q), \ \ k = \max(k_1, \dots, k_{d+1}), \ \ l = \max(l_0,\dots,  l_{d+1}).$$
Below in $3.C$,  we induce extra numbers $M=M_f$ and $c=c_f$ from $f$.
Then the function $f$  and the approximation data determine the constants:
$$(M,c,L,k, D).$$
We say that 
$P(u,u_x,  u_s, u_{xs}, \dots , u_{\alpha s}) $ is
 {\em approximable  in the class} $(M,c,L,k, D)$ in ${\bf PDE}_2^*$,
  if it is induced from some discrete dynamics as above whose constants 
  are all less than or equal to $M,c,L,k,D$ respectively.

 \vspace{3mm}
 
\begin{rem} Notice that three data, defining equations of discrete dynamics 
 (1) $z_{N+1}^{t+1}= f(z_{N-l_0}^{t+1}, \dots,  z_{N+k_{d+1}}^{t-d})$, (2)
 the exponents of the scaling change of variables $(m,p,q)$ and (3) the order to take 
 the Taylor expansions $\alpha$, 
 determine the defining PDEs.
 \end{rem}

\vspace{3mm}

{\bf 3.C Scale transform and tropical geometry:}
A relative $(\max, +)$-function $\varphi$ 
is a  piecewise linear function of the form: 
$$\varphi(\bar{x})=   
 \max(\alpha_1 + \bar{a}_1 \bar{x}, \dots , \alpha_m+ \bar{a}_m \bar{x})
 - \max(\beta_1 + \bar{b}_1 \bar{x}, \dots , \beta_l+ \bar{b}_l \bar{x})$$
 where
 $\bar{a}_l \bar{x}= \Sigma_{i=1}^n a_l^i x_i$, 
 $\bar{x} =(x_1, \dots,x_n) \in {\mathbb R}^n$,  
$ \bar{a}_l =(a_l^1, \dots, a_l^n) , \bar{b} \in {\mathbb Z}^n$
 and $\alpha_i , \beta_i  \in {\mathbb R}$. 
 Notice that $\varphi$ is Lipschitz, since it is piecewise linear.
We  say that the multiple integer $M \equiv ml$
 is the {\em number of the components} of $\varphi$.

Correspondingly
tropical geometry associates  the parametrized  rational function given by:
$$f_t(\bar{z}) \equiv \frac{k_t(\bar{z})}{h_t(\bar{z})}=
\frac{\Sigma_{k=1}^m t^{\alpha_k} \bar{z}^{\bar{a}_k}}
{\Sigma_{k=1}^l t^{\beta_k} \bar{z}^{\bar{b}_k}}$$
 where
$  \bar{z}^{\bar{a}}= \Pi_{i=1}^n z_i^{a^i}$, 
$\bar{z}=(z_1, \dots, z_n) \in {\mathbb R}^n_{>0}$.
We say that $f_t$ above is a {\em relatively elementary} function.
We say that both terms $h_t(\bar{z})= \Sigma_{k=1}^l t^{\beta_k} \bar{z}^{\bar{b}_k}$ 
and $k_t(\bar{z})= \Sigma_{k=1}^m t^{\alpha_k} \bar{z}^{\bar{a}_k}$ are
just elementary functions.

These two functions $\varphi$ and $f_t$ admit one to one correspondence between their presentations.
Moreover the defining equations are transformed  by two steps, 
firstly by  taking conjugates by $\log_t$
and secondly by letting $t \to \infty$  (see [Mi]).

In some cases the same $(\max,+)$ function admits 
different  presentations, 
while the corresponding rational functions are mutually different.
For example for $\varphi(x)\equiv  \max(0,2x) =   \max(0,x,2x)\equiv  \psi (x)$, 
  the corresponding rational functions $f_t(z)= z^2+1$ and $g_t(z)=z^2+z+1$ 
are mutually different.
\begin{df}
Let $f_t$ and $g_t$ be the relatively elementary functions 
with respect to $\varphi$ and $\psi$ respectively.

We say that $f_t$ and $g_t$ are mutually tropically equivalent, 
if $\varphi$ and $\psi$ are the same as maps (but possibly with the different presentations).
\end{df}

For a relatively elementary function $f_t$,
let $c_f \geq 1$ be the Lipschitz constant 
and $M_f$ be the number of the  components 
with respect to the corresponding $(\max, +)$-function $\varphi$.

\vspace{3mm}

{\bf 3.D Asymptotic estimates between different PDEs:}
Let $f$ be a relatively elementary function of $n$ variables.
We have two constants $M_f,c_f \geq 1$ in $3.C$.
 Let 
 $\{(z_{N-l_0}^{t+1}, \dots, z_N^{t+1}, z_{N- l_1}^t ,\dots, z_{N+k_1}^t, \dots, z_{N+k_{d+1}}^{t-d}) ,
(m,p,q)\}$ be an approximation data, and $L,k,  D$ be the corresponding set of the numbers.
With these numbers, consider:
 $$M=\max(M_f,M_g), \quad c = \max(c_f,c_g), \quad C $$
 where $C$ is the maximum of the  error constants for the approximations of PDEs ($3.B$).
We have obtained  the   following result:

\begin{thm}[K2]
Let $f$ and $g$ be both  relatively elementary and increasing functions
of $n $ variables, which are mutually  
tropically equivalent. Let $P, Q \in {\bf PDE}_2^*$ be two PDEs of order 
 $\alpha  \geq 0$ which are approximated by $f$ and $g$
with the above approximation data respectively,
where $*=$ fin or $\infty$.

Let us take positive $C^{\alpha +1}$ solutions $u, v: (0, A_0) \times [0, T_0) \to (0, \infty)$ with: 
$$P(\epsilon, u,u_x, u_s, \dots, u_{\alpha x}, u_{\alpha s})=0, \
Q(\epsilon, v,v_x, v_s, \dots, v_{\alpha x}, v_{\alpha s})=0$$
and assume the estimates  $K_{\alpha+1}(u), K_{\alpha+1}(v) \leq K$.
Then  for any $0 < \epsilon \leq  \min( \frac{1}{2CK}, (L+1)^{-1}A_0, (L+1)^{-1}T_0, n^{-1})$,
they satisfy the asymptotic estimates:
$$(\frac{u(x,s)}{v(x,s)})^{\pm 1}
 \leq  
 (2 M)^{8  \frac{c^{ \epsilon^{- D} (x+ks)   +1}  -1}{c-1}}  
([u:v]_{(L+1)\epsilon})^{c^{ \epsilon^{-D}(x+ks) +n}}.$$

If both $f$ and $g$ happen to be linear, then $c=1$ and so they admit the exponential asymptotic estimates.
\end{thm}
The core of such asymptotic estimates has also appeared in [K1] in the case of  discrete dynamics.

\begin{cor}
Under the above situation, $P_1 \sim_*^+ P_2$  in $(M',c', D,L+1; \alpha)$,
 where:
 $$(M',c') =
 \begin{cases}
 ((2M)^{\frac{8}{c-1}}, c^k) & c>1, \\
 ((2M)^{8k} , 1) & c=1
 \end{cases}$$
 \end{cor}
{\em Proof:}
Firstly suppose $c>1$. 
Then one has the estimates:
\begin{align*}
&  (2 M)^{8  \frac{c^{ \epsilon^{- D} (x+ks)   +1}  -1}{c-1}}  \leq 
 [(2M)^{\frac{8}{c-1}}]^{( c^k)^{\epsilon^{-D}(x+s+1)}}, \\
 & c^{ \epsilon^{-D}(x+ks) +n } \leq (c^k)^{ \epsilon^{-D} (x+s+1)}.
 \end{align*}

 Next suppose $c=1$. Then: 
 $$\lim_{c \to 1}  (2 M)^{8  \frac{c^{ \epsilon^{- D} (x+ks)   +1}  -1}{c-1}} 
 = (2M)^{8 ( \epsilon^{- D} (x+ks)   +1)} \leq [(2M)^{8k}]^{ \epsilon^{- D} (x+s+1)   }.$$
 This completes the proof.

\vspace{3mm}

{\bf 3.D.2 Example:}
Let us apply the above method to a concrete case.
See [K2] for more examples of order $2$.

\begin{prop}  
Two partial differential equations of order $1$:
$$v_s+ \frac{\epsilon}{2} vv_x - \frac{1}{2} v^2=0, \quad 2 u_s  + \frac{\epsilon}{2} u(u_s+ u_x)  =0$$
are both in ${\bf PDE}_2^{\infty}$ in the class $(20,2,1,2,1)$,
and they are mutually related in $\sim_{\infty}^{e^e}$ in the class  $(40, 4, 1,2; 1)$
\end{prop}
{\em Proof:}
We have verified that  these equations are in ${\bf PDE}_2^{\infty}$ in [K2].
For convenience, we give the explicit approximation data here.
Let us choose:
\begin{align*}
& z_{N+1}^{t+1}  = f(z_{N-1}^{t+1} ,z_N^t, z_{N+2}^t) \equiv
\frac{z_{N+2}^t}{2} + 
\frac{ z^t_N(1+2 z_{N-1}^{t+1})}{2(1+ z_N^t)}, \\
&  w_{N+1}^{t+1}  = g(w_{N-1}^{t+1}, w_N^t, w_{N+2}^t) \equiv
\frac{w_{N+2}^t}{2} + 
\frac{ w^t_N+  w_N^t w_{N-1}^{t+1}}{2(1+w_N^t)},
\end{align*}
and $(m,p,q)=(1,1,1)$. Then it is direct to check the followings:
\begin{align*}
& \epsilon a(x+  \epsilon,s+ \epsilon)-   A(\epsilon a(x-  \epsilon,s+\epsilon), \epsilon a(x,s),
\epsilon a(x+2 \epsilon , s) ) \\
& = 
 \begin{cases}
  \frac{\epsilon^2}{2(1+   \epsilon v) }
  (2v_s +2 \epsilon vv_x - v^2) +\epsilon^2 {\bf F}^1(v,v_x, ..,v_{xx}) & (a,A)= (v,f), \\
   \frac{\epsilon^2}{2(1+   \epsilon u) }(2u_s + \epsilon uu_s + \epsilon  uu_x)
+ \epsilon^2 {\bf G}^1(u,u_x, .. ,u_{xx}) & (a,A)=(u,g),
\end{cases}
\end{align*}
where the estimates hold:
\begin{align*}
& |{\bf F}^1|  \leq   \frac{2(\epsilon v)
 \epsilon  |   v_{xs}(\xi_{-11}) | }{(1+   \epsilon v) }
+\epsilon (2 |v_{xs}(\xi_{11}) |+| v_{xx}(\xi_{20})| ) 
   \leq 5 \epsilon ||v||_2 , \\
& |{\bf G}^1|
  \leq   \frac{\epsilon u
 \epsilon  |  u_{xs}(\eta_{-11}) | }{(1+   \epsilon u) }
+\epsilon ( 2|u_{xs}(\eta_{11}) |+|  u_{xx}(\eta_{20})| ) 
   \leq 4 \epsilon ||u||_2 .
   \end{align*}
Thus both the equations are in ${\bf PDE}_2^{\infty}$, and their error constants are both bounded by $5 \epsilon \leq 5$.
Thus combining with theorem $3.1$, one has verified:
\begin{lem}
Let us put $K = \max(K_2(u), K_2(v))$ and
choose any $0 < \epsilon \leq \min(0.1 K^{-1}, 1/3)$. Then any 
 $C^2$ positive solutions to the above two PDEs 
 satisfy the asymptotic estimates for all  $(x,s) \in  (0, \infty)  \times [0, T_0) $: 
$$(\frac{u(x,s)}{v(x,s)})^{\pm 1}
 \leq  
 40^{ 2^{\epsilon^{-1} (x+2s) +4 }  }  
([u:v]_{2\epsilon})^{2^{ \epsilon^{-1}(x+2s) +3}}.$$
 \end{lem}
 
 {\em Proof of proposition $3.4$:}
 The above estimates are bounded by:
$$
 40^{ 2^{\epsilon^{-1} (x+2s) +4 }  }  
([u:v]_{2\epsilon})^{2^{ \epsilon^{-1}(x+2s) +3}}   \leq
 40^{ 2^{2 \epsilon^{-1} (x+s+1)  }  }  
([u:v]_{2\epsilon})^{2^{ 2\epsilon^{-1}(x+s+1) }}.$$

 for all $0 < \epsilon \leq \min(0.1 K^{-1}, 1/3)$.
 This completes the proof.

\section{Unrelated classes}
 Let us take the discrete dynamics $z_{N+1}^{t+1}= f(z_{N-l_0}^{t+1}, \dots,  z_{N+k_{d+1}}^{t-d})$
  and $(m,p,q)$. Then 
 one obtains a family of PDEs $\{ P_{\alpha}\}_{\alpha \geq 1}$ with respect to the order of the Taylor expansions
 (Remark $3.1$ in $3.B$).
 
 Let us fix $(m,p,q)$, and take two discrete dynamics
 $z_{N+1}^{t+1}= f(z_{N-l_0}^{t+1}, \dots,  z_{N+k_{d+1}}^{t-d})$ and
 $w_{N+1}^{t+1}= g(z_{N-l_0}^{t+1}, \dots,  z_{N+k_{d+1}}^{t-d})$.
Then one obtains two families of PDEs $\{P_{\alpha}\}_{\alpha \geq 1}$ 
and $\{Q_{\alpha}\}_{\alpha \geq 1}$ correspondingly.

\begin{df}
Two discrete dynamics 
$z_{N+1}^{t+1}= f(z_{N-l_0}^{t+1}, \dots,  z_{N+k_{d+1}}^{t-d})$ and
 $w_{N+1}^{t+1}= g(z_{N-l_0}^{t+1}, \dots,  z_{N+k_{d+1}}^{t-d})$
 are infinitely unrelated   in $\sim_*^+$, if for any constants 
 $M,c,D,L$, there are some $\alpha$ so that
 $P_{\alpha}$ and $Q_{\alpha}$ are mutually unrelated
 in $\sim_*^+$ in the  class $(M,c,D,L; \alpha)$.
 \end{df}

Here we give some pairs of PDEs which are mutually unrelated.
We treat two cases, where:

\

(1) they are both in ${\bf PDE}_2^{fin}$, which are exponentially unrelated  in $\sim_{fin}^e$.
They arose from two discrete dynamics which are mutually infinitely unrelated.

(2) they are both in ${\bf PDE}_2^{\infty}$, which are double-exponentially unrelated  in $\sim_{fin}^{e^e}$.

\vspace{3mm}

{\bf 4.A Unrelated pairs in the linear case:}
Below we show that the following linear equations are mutually unrelated:

      \begin{thm}   For any $M,D,L >1$, 
  there are  $l_0$ and $a_0,b_0$ so that for each  even $l=2m  \geq l_0$,
  two PDEs $u,v: (0, a_0) \times [0, b_0) \to (0, \infty)$ of order $l$:
  \begin{align*}
&  \epsilon u_s + \frac{\epsilon^2}{2} u_{2s} + \dots + \frac{\epsilon^l}{l!} u_{ls} + \epsilon^2 u_{xs} 
+\frac{\epsilon^3}{2} u_{xss} + \dots + \frac{\epsilon^l}{l!} u_{(l-1)xs}=0, \\
&  \epsilon(v_s+v_x) + \epsilon^2 (v_{xs} + \frac{1}{2}v_{2s} + \frac{1}{2} v_{2x}) + \dots + \frac{\epsilon^l}{l!} v_{lx} =0
\end{align*}
 satisfy the followings;
 
 (1) they are approximable  in the class $(1, 1, 0, 1, 1)$ in ${\bf PDE}_2^{fin}$, and
 
 (2) they 
  are mutually unrelated in   $(M,1, D, L; l)$.
 \end{thm}
  {\em Proof:}
  (1) Let us consider the discrete dynamics  given by:
  $$z_{N+1}^{t+1}=z_{N+1}^t, \quad w_{N+1}^{t+1}=w_N^t$$
   and the scaling parameters  by
   $N = \frac{x}{\epsilon},  t = \frac{s}{\epsilon}$ and $ z_N^t = u(x,s)$.
  Then by taking the Taylor expansions up to order $l+1$, one obtains
  the desired PDEs, which are both approximable in ${\bf PDE}_2^{fin}$
  in the class $(1,1,0, 1,1,1)$.
 Notice that each monomial of the  PDE in $u$  contains derivatives of $s$.

(2) 
  Let $f(x) = x^l+1$, and put $u(x,s) = f(x)$ and $v(x,s) = f(x-s)$.
  It is immediate to see that they are solutions respectively, 
  because of independence of the variable $s$ for $u$, and
  of the symmetry of the equation for $v$.
  Moreover   the equalities $K(u)_{l+1}=K_{l+1}(v)=0$ hold,
  since $l+1$ derivatives of $u$ and $v$ are both equal to zero.
  We require that $l$ are even in order to guarantee positivity of values of $v$.

  If they were related in  $(M,1,D, L; l)$, 
  then there is some $C$ independent of solutions 
  so that  they must satisfy the asymptotic estimates:
  $$(\frac{u(x,s)}{v(x,s)})^{\pm 1}
 \leq   M^{ \epsilon^{- D} (x+s+1) }  [u:v]_{L \epsilon}$$
  for all $0 < L \epsilon \leq \min(a_0, b_0,C)$.

  Let us choose small $\epsilon >0$ with $L \epsilon \leq 1$. Then the estimates hold:
  \begin{align*}
&  M^{ \epsilon^{- D} (x+s+1)} \leq  M^{ \epsilon^{- D} (a_0+b_0   +1)} , \\
&   [u:v]_{L \epsilon} \leq [u:v]_1= \max(\sup_{0 < x \leq a_0}  \frac{x^l+1}{(x-1)^l +1}, \ 
\sup_{0 < x \leq 1} \frac{(b_0-x)^l+1}{x^l +1})
  \end{align*}
  where:
  \begin{align*}
&    \frac{x^l+1}{(x-1)^l +1} = \frac{1+ \frac{1}{x^l}}{(1- \frac{1}{x})^l + \frac{1}{x^l}}  \leq
    \begin{cases} \frac{2}{\frac{1}{2^l}} = 2^{l+1} & x \geq 2, \\
    2^{l+1} & x \leq 2.
    \end{cases}, \\
    & \sup_{0 < x \leq 1} \frac{(b_0-x)^l+1}{x^l +1} = b_0^l +1.
    \end{align*}
    Let us choose $b_0 \geq 3$. Then 
     the estimate 
    $ [u:v]_{L \epsilon} \leq b_0^l +1$ holds.

    Now let us find $a_0 >b_0 \geq 3$ so that
    the inequality 
    $\frac{a_0}{a_0-b_0} > b_0$ holds.
    In fact  for $2< \beta < 4$ (say $\beta =3$ is enough),
    if one chooses large $a_0 $ with $a_0^2 - \beta a_0 \geq 0$
    then 
    $$b_0 = \frac{1}{2}(a_0 + \sqrt{a_0^2 - \beta a_0}) \geq 3$$
    satisfies the required conditions. 
        Notice that  the equality $\frac{a}{a-b}=b$ holds,
         in the case when $b=  \frac{1}{2}(a + \sqrt{a^2 - 4 a})$.
    Thus if one chooses $b_0 $ as above, then the estimate     $\frac{a_0}{a_0-b_0} > b_0$
     holds    by an elementary observation.

    Let us fix such a pair $(a_0,b_0)$.
   Now  $\frac{u(a_0,b_0)}{v(a_0,b_0)}= \frac{a_0^l +1}{(a_0-b_0)^l +1}  $ holds.
    If one chooses sufficiently large $l >>1$,
  then the estimates: 
  \begin{align*}
    \frac{u(a_0,b_0)}{v(a_0,b_0)} (b_0^l+1)^{-1} 
 & =
  \frac{a_0^l +1}{(a_0-b_0)^l +1}\frac{1}{ (b_0^l+1)}  \ \
   \begin{cases}
  \geq \frac{1}{4} \frac{a_0^l }{b_0^l} & (a_0 -b_0 \leq 1) \\
  \geq \frac{1}{4} \frac{a_0^l}{(a_0-b_0)^l} \frac{1}{b_0^l} & (a_0-b_0 >1)
  \end{cases} \\
&   >   M^{\epsilon^{- D} (a_0+b_0   +1)}
  \end{align*}
    hold, since the last term is independent of $l$.
   Then  one has:
    $$ (\frac{u(a_0,b_0)}{v(a_0,b_0)})^{\pm 1}  >
      M^{\epsilon^{- D} (a_0+b_0 +1)}(b_0^l +1)
      \geq    M^{ \epsilon^{- D} (a_0+b_0   +1)}[u:v]_{L \epsilon}.$$
      This is a contradiction.
  This completes the proof.

  \vspace{3mm}

 {\bf 4.B Non linear estimates:}
 Let us treat the case of double-exponential estimates.
 Let $P(\epsilon, u,u_s,u_x, \dots, u_{\alpha s})=0$ be in ${\bf PDE}_2^{\infty}$
 of order $\alpha \geq 1$,
 and let us compare its solutions with the translations: 
 $$P_I: v_s =Iv_x \qquad (I >0).$$

Let us start from a general situation.
 
\begin{lem}
 Suppose that for some $\frac{1}{4} \geq \delta_0>0$, $C,C' , C''> 0$ 
 and  for all small $0< \epsilon \leq \epsilon_0$,
   there are solutions $u_{\epsilon}(x,s)$ on $(0, \delta_0) \times [0, \delta_0) \to (0, \infty)$
   with the initial values $u(x,0)=f(x)=1-x$, which satisfy both the estimates:
   $$  (1) \ \  C \leq u_{\epsilon}(x,s) \leq C', \quad  (2)  \ \ K_{\alpha+1}(u_{\epsilon}) \leq C''.$$
 Then for any $M,c,D,L \geq1$,  there is some $I_0 >0$ so that for all $I \geq I_0$,
 $P$ and $P_I$ are unrelated in $\sim_{fin}^{e^e}$ in the class  $(M,c,D,L; \alpha)$
 \end{lem}
 {\em Proof:}
 We verify the conclusion for a specific $I_0>0$. Then the general case follows by restricting
 the domain of the solutions by choosing some smaller $\delta_0 \geq \delta_0' >0$.
 
 Let  us put $v(x,s)= f(x+I_0s)$. Then $v$ is the solution to $P_{I_0}$ and $K_{\alpha+1}(v)=0$ holds.
 
 Let us take  sufficiently small $\delta_0 >> \epsilon >> \delta>0$, and choose $I_0$ with
 $I_0 \delta_0 = 1- \delta_0 -\delta$.
 One may assume  the estimate
  $I_0 L \epsilon \leq \delta_0$. Then
  $\delta_0 + I_0 L \epsilon \leq 2 \delta_0$ and 
  $L \epsilon + I_0 \delta_0 = 1- \delta_0 - \delta + L \epsilon \leq 1 - \frac{ \delta_0}{2}$. So
  the estimate
  $x+ I_0 s \leq 1-  \frac{\delta_0}{2} $
   holds on the initial domain
$ (x,s) \in (0, \delta_0) \times [0, L\epsilon) \cup (0, L\epsilon) \times [0, \delta_0)$.
 Thus the estimate holds:
  $$[v:u_{\epsilon}]_{L \epsilon} \leq 2C' \delta_0^{-1}.$$

Suppose they could be  $(M, c, D,L;\alpha)$ related.
Then they must satisfy the asymptotic estimates
$(\frac{u_{\epsilon}(x,s)}{v(x,s)})^{\pm 1}
 \leq   M^{ c^{2 \epsilon^{- D} } }  (2C' \delta_0)^{c^{ 2\epsilon^{-D}}}$.

On the other hand $v(\delta_0,\delta_0) = \delta$ and $u_{\epsilon}(\delta_0, \delta_0) \geq C$ hold.
So $\frac{u_{\epsilon}(\delta_0, \delta_0)}{v(\delta_0, \delta_0)} \geq \frac{C}{\delta}$ which can be
arbitrarily large. This is a contradiction.
This completes the proof.

\vspace{3mm}

{\bf 4.B.2 Conservation equations:}
Let us apply the above situation to the first order conservation equations and
 the translations:
 $$ 4u_s +\epsilon u(u_s+u_x)=0 \ , \quad  v_s=Iv_x  .$$ 
We verify that both are in ${\bf PDE}_2^{\infty}$ and are mutually unrelated in $\sim_{fin}^{e^e}$
(Compare this with proposition $3.3$):
\begin{thm}
For any $M, c, D, L$, there is $I_0$ so that for all $I \geq I_0$,
 pairs of PDEs given  by:
$$ 4u_s +\epsilon u(u_s+u_x)=0,  \quad v_s=I v_x$$
are both in ${\bf PDE}_2^{\infty}$  in $(20, 2,1,I+1, 1)$, and 
are unrelated in $\sim_{fin}^{e^e}$ in the class  $(M, c, D, L;1)$.
\end{thm}
Before proceeding, let us briefly recall a way to produce solutions,
  called  the {\bf method of characteristics} for 
   the conservative non linear equations
 of the form $u_s + F(u)u_x =0$.
 
 Let $u(\xi,0)=f(\xi)$   be the initial condition, and
 try to solve the equation $\frac{dx}{ds}= F(u)$ for $u=u(x,s)$ with $x(\xi, 0)= \xi$.
 Then $\frac{d u(x,s)}{ds}=u_s + \frac{dx}{ds}u_x =0$ holds. So $u$ is constant along $x(\xi, s)$.
  Moreover $\frac{d^2 x}{ds^2}= \frac{F(u)}{ds} \frac{du}{ds}=0$, and so 
  $x(\xi,s) = F(f(\xi))s + \xi$. 
 Thus   if one could solve $\xi =\xi(s,x)$, then $f(\xi(x,s))$ will give us solutions, 
since $u$ is constant  along $x(\xi,s)$, 
\vspace{3mm}

{\em Proof:}
For the second equation, it is induced from the discrete dynamics
$w_{N+1}^{t+1}= w^t_{N+I+1}$, so it lies in the class
$(1, 1,0, I+1, 1)$.
Thus both are in ${\bf PDE}_2^{\infty}$  in the class $(20, 2,1,I+1, 1)$ by proposition $3.3$.

Now let us consider the initial function $f: (0,1) \to (0, \infty)$  by  $f(x)= 1-x$.
The result follows if one can find solutions $u$ which satisfy two conditions (1) and (2) in lemma $4.2$.
We will solve the equation very concretely by use of the method of characteristics.

Let us consider the equation:
$$x-\xi =\frac{\mu f(\xi)s}{1+ \mu f(\xi)} = \frac{\mu(1- \xi)s}{1+ \mu(1-\xi)}, \qquad 
( \ \mu = \frac{\epsilon}{4} \ ).$$
Then it gives the equation  \
$\mu \xi^2 -(1+ \mu(x-s+1))\xi +(1+\mu)x-\mu s=0$, and one can solve it with $x(\xi, 0)=\xi$:
$$\xi (x,s)= \frac{1}{2\mu} \ [ \ 1+\mu(x-s+1)- \sqrt{(1+\mu(x-s+1))^2 - 4\mu \{(1+\mu)x- \mu s\} } \ ].$$

Let us put the solution $u(x,s) = 1- \xi (x,s)$.
Then for a small $\delta_0 >0$ and all sufficiently small $\epsilon >0$, 
 two conditions 
  (1)  $C \leq u(x,s) \leq C'$ and  (2)  $K_2(u) \leq C''$ are certainly satisfied,
   by elementary calculations.
   
This completes the proof.

\vspace{3mm}

\begin{rem}
One has an observation of positivity of solutions, which rely heavily on 
the structure of the equation.
Let us rewrite the equation as $u_s = -     \epsilon u u_x(4+\epsilon u)^{-1}$.
At $s=0$, $u_x(x,0) =-1 <0$ hold and so 
$u_s(x,0) = \epsilon u (4+ \epsilon u)^{-1} >0$ holds.
Thus there is some $T_0  >0$ so that
$u_s(x,s) >0$ still hold for all $(x,s) \in (0,\frac{1}{2}) \times [0, T_0)$.
In particular $u(x, s) \geq \frac{1}{2}$ for all $(x,s)  \in (0, \frac{1}{2}) \times [0, T_0)$.
The same argument works for 
the equation $v_s+ \frac{\epsilon}{2}vv_x - \frac{1}{2}v^2=0 $.
\end{rem}

\section{Computational aspects}
Our mathematical framework is intimately familiar with the computer systems.
Many essential quantities of the analysis
are computable within finite steps and finite values.
By use of numerical calculations, 
we reprove  Theorem 1.3 for some particular values of constants
($M=10^3$ and $l =10^3$).
This can be done possible since the discretization of the PDEs are the rigorous 
solutions rather than just approximating. However in more general situations,
such verification would be expected to work effectively by use of 
numerical simulation with guaranteed accuracy.

For the proof that the two partial differential equations are of
unrelated pair, it is sufficient to 
show that there exists the test point does not satisfy the inequality
 of the asymptotic estimates in the definition 2.1.
We focus on this inequality and give some discussions from the
point of view of the numerical calculations in two cases:
the related pair and the unrelated pair.

 Recall that we have induced PDEs from discrete dynamics, 
 and our computation here will be done for these discrete dynamics.

 Let $\{z_N^t\}_{N,t \geq 0} $ and  $\{w_N^t\}_{N,t \geq 0} $ be two discrete dynamics.
 Then for $L, N_0$ and $t_0$, we put the discrete version of the initial rates by:
 $$[\{z_N^t\} : \{w_N^t\}]_{L, N_0,t_0} \equiv 
 \sup_{0 \leq N \leq N_0, 0 \leq t \leq t_0, 0 \leq a \leq L} 
 \max \{ (\frac{z_N^a}{w_N^a})^{\pm} ,  (\frac{z_a^t}{w_a^t})^{\pm}\}
 $$
 Let us denote $\tilde{u}(\epsilon N, \epsilon t) =z_N^t$ and
  $\tilde{v}(\epsilon N, \epsilon t) =w_N^t$ respectively.
Then we regard that both $\tilde{u}$ and $\tilde{v}$ approximate $u$ and $v$ respectively:
$$\tilde{u}(\epsilon N, \epsilon t) \sim u(\epsilon N, \epsilon t) $$
and similar for $v$. In fact in the examples in $5.A$,  
both $\tilde{u}$ and $\tilde{v}$ coincide with $u$ and $v$ respectively, which induces 
Proposition $5.1$ below.

Throughout this section, we choose the rescaling parameters 
$(m,p,q)=(1,1,1)$.
All of the numerical calculations in this section are performed by
using the computer algebra system ``Maple 13'' with rational or floating
number
manipulations.
Then numbers after calculation are converted to the floating-point
numbers for presentation purposes.

\vspace{3mm}

{\bf 5.A Case of the unrelated pair:}
Let us recall two PDEs:
  \begin{align*}
&  \epsilon u_s + \frac{\epsilon^2}{2} u_{2s} + \dots + \frac{\epsilon^l}{l!} u_{ls} + \epsilon^2 u_{xs} 
+\frac{\epsilon^3}{2} u_{xss} + \dots + \frac{\epsilon^l}{l!} u_{(l-1)xs}=0,  \qquad (*_1)\\
&  \epsilon(v_s+v_x) + \epsilon^2 (v_{xs} + \frac{1}{2}v_{2s} 
+ \frac{1}{2} v_{2x}) + \dots + \frac{\epsilon^l}{l!} v_{lx} =0 \qquad (*_2)
\end{align*}
 which are mutually induced by the discrete dynamics:
 $$z_{N+1}^{t+1}=z_{N+1}^t \ , \ w_{N+1}^{t+1}=w_N^t.$$
 
 Recall that in $4.A$ we have verified that the pair of the solutions 
 $u(x,s)=x^l +1$ and $v(x,s)=(x-s)^l +1$ breaks the exponential bounds:
   $$\max\{ \ \frac{u(x,s)}{v(x,s)} \ ,  \  \frac{v(x,s)}{u(x,s)} \ \}
 >>   M^{ \epsilon^{- D} (x+s+1) }  [u:v]_{L \epsilon}.$$

Now  for $N,t = 0,1,2, \dots$, 
let us put the solutions to the above discrete dynamics:
$$z_N^{t}=(\epsilon N)^l+1 \ , \  w_N^{t}=\epsilon^l(N-t)^l+1$$
 which are precisely the same as the solutions to the PDEs $u(x,s)$ and $v(x,s)$ 
 by  the scaling change of variables  $x=N\epsilon$ and $s= t \epsilon$
 respectively.

 We put:
 $$\tilde{Q}_e(x,s) =  \log \max\{ \frac{w_N^t}{z_N^t} , \frac{z_N^t}{w_N^t} \} 
  - \log \ [\{z_N^t\} : \{w_N^t\}]_{L, \epsilon^{-1}x, \epsilon^{-1}s}$$
 whose values are much bigger than $   \epsilon^{- D} (x+s+1) \log M$ in  the above case.

We calculate the values $Q_e(x,s)$ at the points
$(x,s)=(N\epsilon, t\epsilon)$ for $\{N,t \in {\mathbb Z}\,|\, 1< N \le
A_0/\epsilon, 1<t \le T_0/\epsilon\}$.

Let us put:
 \begin{align*}
Q_e(x,s)= \max\left(0,  \tilde{Q}_e(x,s) \right).
\end{align*}

Now we verify the following by computer calculations.
These calculations are performed by using rational numbers which 
means that the calculated values are exact ones.
 \begin{prop}{\label{prop:comp}}
 The pair of the above equations $(*_1,*_2)$ are mutually unrelated 
 in $\sim_{fin}^e$ in $(10^3,1,1,1; 10^3)$ at $\epsilon =1/2$.
 \end{prop}
 {\em Proof:}
 Since the values of discrete dynamics $z_N^t$ and $w_N^t$ are 
 precisely the same as the solutions $z_N^t = u(N\epsilon, t\epsilon)$
 and $w_N^t = (x,s)=(N\epsilon, t\epsilon)$ respectively,
  it is enough to verify that 
 $Q_e(x,s)$ certainly hit bigger values than   $   \epsilon^{- D} (x+s+1) \log M = 6(x+s+1) \log 10$
 at some points.

Let us choose constants by $l=1000,  L=0 $ and $\epsilon=1/2$, and
consider the values of $\tilde{Q}_e$ 
at each point $(x,s)=(N\epsilon,t\epsilon)$ 
for $1 \leq N,t \leq 8$.
For example 
$\tilde{Q}_e(1.5,2) = 
\log \max\{ \frac{w_3^4}{z_3^4} , \frac{z_3^4}{w_3^4} \}  - \log \ [\{z_N^t\} : \{w_N^t\}]_{0, 3, 4}$

Then we show the values of $Q_e(x,s)$ in 4-digits of precision at the point $(x,s)=(N\epsilon,t\epsilon)$ for $\{N,t \in {\mathbb
Z}\,|\, 1 < N,t\le 8= 4/\epsilon\}$
in Table 1.

\begin{table}[htbp]
\caption{Computed values of $Q(N\epsilon,t\epsilon)$ (Unrelated pair)}
\label{table1}
\begin{center}
\begin{tabular}{c|ccccccccc} 
\hline
$N\backslash t$ &0&1&2&3&4&5&6&7&8\\
\hline
0&0.0 &0.0&0.0&0.0&0.0&0.0&0.0&0.0&0.0\\
1&0.0&0.0&0.0&0.0&0.0&0.0&0.0&0.0&0.0\\
2&0.0&0.6931&0.0&0.0&0.0&0.0&0.0&0.0&0.0\\
3&0.0&404.8& 404.8&0.0&0.0&0.0&0.0&0.0&0.0\\
4&0.0&287.7& 691.8& 287.6&0.0&0.0&0.0&0.0&0.0\\
5&0.0&223.1& 510.1& 510.1& 223.2&0.0&0.0&0.0&0.0\\
6&0.0&182.3& 404.8& 287.6& 404.9& 182.7&0.0&0.0&0.0\\
7&0.0&154.2& 335.8& 154.1& 154.2& 335.7& 154.0&0.0&0.0\\
8&0.0&135.5& 287.0& 64.50&0.0& 64.50& 287.0& 133.0&0.0\\
\hline
\end{tabular}
\end{center}
\end{table}

Now for $D=1,M=1000$,
$\log M^{\epsilon^{-D}(x+s+1)}\le\log M^{\epsilon^{-D}(4+4+1)} \sim  124.3$.
So one can find that 
the exponential asymptotic estimates:
\begin{align*}
 Q_e(x,s) \le \log M^{\epsilon^{-D}(x+s+1)}
\end{align*}
 does not hold at the several points.
 This completes the proof.

\vspace{3mm}

\begin{rem}Owing to the exact solutions of the partial differential equations,
we are able to prove {\rm proposition \ref{prop:comp}} by using computer system.
Nevertheless, in any case, even that we know little of solutions,
it might be possible to apply the self-validating numerical method {\rm [Loh,O]}
in mathematical proof.
\end{rem}

\vspace{3mm}

{\bf 5.B Case of the related pair:}
In case of the examples of the related pairs we present here, 
the numerical calculations do not directly provide the mathematical
proof of the relevancy.
Nevertheless such calculations give us 
several insights on  the actual behaviour for further analysis by 
comparison with the case of the unrelated pair.

Let us recall two PDEs:
$$v_s+ \frac{\epsilon}{2} vv_x - \frac{1}{2} v^2=0, \quad 2 u_s  + \frac{\epsilon}{2} u(u_s+ u_x)  =0$$
which are mutually  induced from the discrete dynamics:
\begin{align*}
 & z_{N+1}^{t+1}=\dfrac{z_{N+2}^t}{2} +
   \dfrac{z_{N}^t(1+2z_{N-1}^{t+1})}{2(1+z_{N}^t)},\quad (**_1) \\
&  w_{N+1}^{t+1}=\dfrac{w_{N+2}^t}{2} +
   \dfrac{w_{N}^t(1+w_{N-1}^{t+1})}{2(1+w_{N}^t)} \quad (**_2)
\end{align*}

In Proposition $1.2$, we have seen the asymptotic estimates for all solutions, which are
 equivalent to:
\begin{eqnarray*}
 \log \left(\dfrac{u(x,s)}{v(x,s)}\right)^{\pm 1}
 - 
{2^{\epsilon^{-1}(x+2s)+3}}
 \log([u:v]_{2\epsilon})
 \le
 {2^{\epsilon^{-1}(x+2s)+4}} \log 40
\end{eqnarray*}

Now
let us consider the solutions to the above discrete dynamics $(**)$ 
 with the initial and boundary values:
$$z_N^{t}=(\epsilon N)^l+1 \ , \  w_N^{t}=\epsilon^l(N-t)^l+1$$
respectively, where $(N,t) =\{0\} \times {\mathbb N} \cup {\mathbb N} \times \{0\}$.

Let $\tilde{u}$ and $\tilde{v}$ be as in $5.A$,
and $u,v: [0,A_0] \times [0,T_0]$ be  solutions 
to the corresponding PDEs respectively.
Then as before we regard that both $\tilde{u}$ and $\tilde{v}$ approximate
$u$ and $v$ respectively  at the points
$(x,s)=(N \epsilon^p, t \epsilon^q)=(N \epsilon, t \epsilon)$
for $\{N,t \in {\mathbb{Z}}\,|\,1<N  \le A_0/\epsilon, 1<t  \le T_0/\epsilon$ \}.

  We put:
 $$\tilde{Q}_{e^e}(x,s) =  \log \max\{ \frac{w_N^t}{z_N^t} , \frac{z_N^t}{w_N^t} \}  -
 {2^{\epsilon^{-1}(x+2s)+3}}
   \log \ [\{z_N^t\} : \{w_N^t\}]_{L, \epsilon^{-1}x, \epsilon^{-1}s}.$$
 Here we calculate the values using floating numbers with 100 digits of precision,
\begin{eqnarray*}
Q_{e^e}(x,s)= \max\left(0,  \tilde{Q}_{e^e}(x,s) \right)
\end{eqnarray*}
where
we choose constants 
$l=1000,  L=1 $, $\epsilon=1/10$ and $A_0=T_0=1$.
In particular the domain of PDEs are 
$\{(x,s)=(N\epsilon,t\epsilon) : 0 \leq N,t \leq 10\}$.
 
 The numerical calculations verify that 
 all the entries are equal to $0$. 
  In particular the estimates: 
 $$Q_{e^e}(x,s) \le {2^{4+(x+2s)/\epsilon}} \log 40$$
  follows, which gives 
 the numerical verification of proposition $3.3$ for particular constants.
Notice that approximately the value is given:
\begin{eqnarray*}
  {2^{4+(x+2s)/\epsilon}} \log 40|_{\epsilon =1/10,x=s=10\epsilon}
\sim 6.337 \times 10^{10}
\end{eqnarray*}

So far we have checked that certainly double exponential estimates hold
for these pairs. Next let us examine whether they might still satisfy the exponential estimates.

Let us consider the values of $Q_e(x,s)= \max\left(0,  \tilde{Q}_{e^e}(x,s) \right)$,
where: 
 $$\tilde{Q}_e(x,s) =  \log \max\{ \frac{w_N^t}{z_N^t} , \frac{z_N^t}{w_N^t} \}  -
   \log \ [\{z_N^t\} : \{w_N^t\}]_{1, \epsilon^{-1}x, \epsilon^{-1}s}$$
 with the same  constants, 
$l=1000,  L=1 $, $\epsilon=1/10$ and $A_0=T_0=1$.
Table 2 gives the result of numerical calculations.


\begin{table}[htbp]
\caption{Computed values of $Q_e(N\epsilon,t\epsilon)$ (Related pair)}
\label{table2}
\begin{center}
\begin{tabular}{c|ccccccccccc}
\hline
$N\backslash t$ &0&1&2&3&4&5&6&7&8&9&10\\
\hline
0&0.0& 0.0& 0.0& 0.0& 0.0& 0.0& 0.0& 0.0& 0.0& 0.0& 0.0\\
1&0.0& 0.0& 0.0& 0.0& 0.0& 0.0& 0.0& 0.0& 0.0& 0.0& 0.0\\
2&0.0& 0.0& 0.095& 0.189& 0.258& 0.331& 0.393& 0.463& 0.524& 0.0& 0.0\\
3&0.0& 0.0& 0.196& 0.327& 0.461& 0.567& 0.682& 0.776& 0.0& 0.0& 0.0\\
4&0.0& 0.0& 0.229& 0.450& 0.618& 0.794& 0.926& 0.0& 0.0& 0.0& 0.0\\
5&0.0& 0.0& 0.318& 0.565& 0.813& 0.991& 0.0& 0.0& 0.0& 0.0& 0.0\\
6&0.0& 0.0& 0.334& 0.670& 0.899& 0.0& 0.0& 0.0& 0.0& 0.0& 0.0\\
7&0.0& 0.0& 0.399& 0.677& 0.0& 0.0& 0.0& 0.0& 0.0& 0.0& 0.0\\
8&0.0& 0.0& 0.294& 0.0& 0.0& 0.0& 0.0& 0.0& 0.0& 0.0& 0.0\\
9&0.0& 0.0& 0.0& 0.0& 0.0& 0.0& 0.0& 0.0& 0.0& 0.0& 0.0\\
10&0.0& 0.0& 0.0& 0.0& 0.0& 0.0& 0.0& 0.0& 0.0& 0.0& 0.0\\
\hline
\end{tabular}
\end{center}
\end{table}

Let us compare their values with:
\begin{align*}
 10 \log 40 =    \epsilon^{-1} \log 40|_{\epsilon =1/10} &  \leq \ 
   {{\epsilon^{-1}(x+s+1)}} \log 40 \\
   &   \leq \ 
    {{\epsilon^{-1}(x+s+1)}} \log 40|_{\epsilon =1/10,x=s=10\epsilon}=
   30 \log 40
\end{align*}
where the  left and right hand  sides are 
 approximately $36.89$ and   $110.7$ respectively. Thus 
the inequality above holds for any points $(N,t)\in [0,10]\times[0,10]$.

Finally  we would like to raise a question:

\vspace{3mm}

\noindent{\bf Question 5.1:} Are
 two PDEs:
$$v_s+ \frac{\epsilon}{2} vv_x - \frac{1}{2} v^2=0, \quad 2 u_s  + \frac{\epsilon}{2} u(u_s+ u_x)  =0$$
exponentially related  in $\sim_{\infty}^e$ or in  $\sim_{fin}^e$ ?


\end{document}